\documentclass[12pt]{amsart}

\usepackage[all]{xy}
\usepackage{amssymb}
\newtheorem{theorem}{Theorem}

\newtheorem{lemma}[theorem]{Lemma}
\newtheorem{prop}[theorem]{Proposition}
\newtheorem{cor}[theorem]{Corollary}

\theoremstyle{remark}

\newtheorem{remark}[theorem]{Remark}

\newtheorem*{problem*}{Problem}
\newtheorem*{remark*}{Remark}
\newtheorem*{convention*}{Convention}
\newtheorem*{notation*}{Notation}
\newtheorem*{examples*}{Examples}
\newtheorem*{example*}{Example}
\newtheorem*{warning*}{Warning}

\def\C{{\mathbb C}}

\def\H{{H}}
\def\L{{\mathcal L}}
\def\K{{\mathcal K}}
\def\I{{\mathcal I}}

\newcommand{\newspan}{\operatorname{span}}
\newcommand{\clsp}{\overline{\operatorname{span}}}

\newcommand{\tr}{\operatorname{tr}}
\newcommand{\rank}{\operatorname{rank}}
\newcommand{\Ind}{\operatorname{Ind}}
\newcommand{\Prim}{\operatorname{Prim}}
\newcommand{\dashind}{\operatorname{\!-Ind}}
\begin{document}

\title[Properties preserved under Morita equivalence]{\boldmath{Properties preserved under\\ Morita equivalence of $C^*$-algebras}}

\author[an Huef]{Astrid an Huef}
\address{School of Mathematics\\
The University of New South Wales\\
NSW 2052\\
Australia}
\email{astrid@unsw.edu.au}

\author[Raeburn]{Iain~Raeburn}
\address{Iain Raeburn, School  of Mathematical and Physical Sciences, University of
Newcastle, NSW 2308, Australia}
\email{iain.raeburn@newcastle.edu.au}

\author[Williams]{Dana P. Williams}
\address{Department of Mathematics\\Dartmouth College\\ Hanover, NH 03755\\USA}
\email{dana.williams@dartmouth.edu}
\address{}

\date{28 November 2005}
\begin{abstract}
We show that  important structural properties of $C^*$-algebras and the multiplicity numbers of representations are preserved under Morita equivalence. 
\end{abstract}
\thanks{This research was supported by the Australian Research Council, the National Science Foundation, the Ed Shapiro Fund at Dartmouth College and the University of New South Wales.}

\maketitle

\section*{Introduction}

Morita equivalence for $C^*$-algebras was introduced by Rieffel in the 1970s, and is now a standard tool in the subject. Saying that two $C^*$-algebras are Morita equivalent is a strong way of saying that ``they have the same representation theory'', and hence one expects representation-theoretic properties of $C^*$-algebras to be preserved by Morita equivalence. Here we aim to provide a brief but comprehensive discussion of this issue, thereby updating and extending previous work of Zettl \cite{zettl1,zettl2}.

Our main new results are that the upper and lower multiplicity numbers of Archbold~\cite{a} and the relative multiplicity numbers of Archbold-Spielberg~\cite{as} are preserved by Morita equivalence: if $A$ and $B$ are Morita equivalent, $\pi\in\hat B$, and $\Ind\pi$ is the corresponding representation of $A$, then the multiplicities of $\pi$ and $\Ind\pi$ coincide (Theorem~\ref{mupi} and its corollaries). We also give a short direct proof that nuclearity is preserved, avoiding previous authors' reliance on Connes' equivalence between nuclearity of $A$ and injectivity of $A^{**}$ (see \cite{zettl1,beer}). We have tried to  use only the basic theory of Morita equivalence, as expounded in \cite[Chapter~3]{tfb}, and we have preferred arguments which do not require separability hypotheses. We have therefore resisted temptations to reduce Morita equivalence to stable isomorphism using the Brown-Green-Rieffel theorem (as in \cite[Theorem~5.55]{tfb}, for example).

We prove in \S\ref{secCCR} that liminarity and related properties are preserved, and that the properties of having continuous trace or bounded trace are preserved. Many of these results were first proved by Zettl using similar arguments \cite{zettl2}, and we have included them here partly to provide a convenient reference in modern notation, and partly because we need the main technical results (Lemma~\ref{lem-green} and its corollaries) in the proof of our main theorem in \S\ref{secRob}. In the last section, we prove that nuclearity is preserved.

\subsection*{Background} We say that two $C^*$-algebras are Morita equivalent to mean that there is an $A$--$B$ imprimitivity bimodule. The basic theory of Morita equivalence was developed by Rieffel, and the last two authors provided a detailed account of his theory in \cite[Chapters~2 and~3]{tfb}, which we use as our main reference.  

If ${}_AX_B$ is an imprimitivity bimodule and $\pi$ is a representation of $B$ on a Hilbert space $H$, we denote by $X\dashind \pi$ or $\Ind\pi$ the induced representation of $A$ on $X\otimes_B H$ characterised by $\Ind\pi(a)(x\otimes h)=(a\cdot x)\otimes h$. The kernel of the representation $\Ind\pi$ depends only on the kernel of $\pi$, so there is a well-defined map $\Ind=X\dashind$ from the set of ideals $\mathcal{I(B)}$ of $B$ to $\mathcal{I}(A)$, which turns out to be an inclusion-preserving bijection with inverse $\widetilde X\dashind$ implemented by the dual bimodule ${}_B\widetilde X_A$ \cite[Theorem~3.22]{tfb}. We refer to this bijection as the Rieffel correspondence associated to the imprimitivity bimodule $X$. If $I$ is an ideal in $B$, then $I$ and $\Ind I$ are canonically Morita equivalent, and so are the quotients $A/\Ind I$ and $B/I$ \cite[Proposition~3.25]{tfb}. The map $\pi\mapsto \Ind\pi$ respects unitary equivalence and irreducibility, and induces a homeomorphism of the spectrum $\hat B$ onto $\hat A$. (It is proved in \cite[Corollary~3.33]{tfb} that the Rieffel correspondence gives a homeomorphism of $\Prim B$ onto $\Prim A$, and it follows from this and the definition of the topology on the spectrum given in \cite[Definition~A.21]{tfb} that it is also a homeomorphism on spectra.)

\section{Properties associated to the algebra of compact operators}\label{secCCR}

A $C^*$-algebra is \emph{elementary} if it is
isomorphic to the algebra $\K(H)$ of compact operators on some Hilbert space $H$. 

\begin{prop}\label{elementary}
Suppose that ${}_AX_B$ is an imprimitivity bimodule. Then $A$ is elementary if and only if $B$ is elementary.
\end{prop}

\begin{proof}
Suppose $B$ is elementary. The algebra $\K(H)$ is Morita equivalent to $\C$ \cite[Examples~2.11 and 2.27]{tfb}, and Morita equivalence is an equivalence relation \cite[Proposition~3.18]{tfb}, so $A$ is Morita equivalent to $\C$; let ${}_AY_{\C}$ be an $A$--$\C$ imprimitivity bimodule. Since a Hilbert $\C$-module $Y_{\C}$ is a Hilbert space and $\K(Y_{\C})$ is then the usual algebra of compact operators \cite[Example~2.27]{tfb}, we deduce that $A=\K(Y_{\C})$ is elementary.
\end{proof}

Following \cite{dix}, we say that a $C^*$-algebra $A$ is \emph{liminary}\footnote{We have followed our sadly missed friend Gert Pedersen in avoiding the dreaded ASHCEFLC (see \cite[\S6.2.13]{ped}), and in preferring to translate the French word \emph{liminaire} as liminary in parallel with the obvious translation of \emph{pr\' eliminaire}.} if $\pi(A)\cong \K(\H_{\pi})$ for all $\pi\in\hat A$, and
\emph{postliminary} if every non-zero quotient of $A$ has a non-zero liminary~ideal. 

\begin{prop}
\label{prop-me-ccr}
Suppose that ${}_AX_B$ is an imprimitivity bimodule.  Then $A$ is
liminary if and only if $B$ is liminary.
\end{prop}

\begin{proof} Assume $A$ is liminary, and let $\pi\in\hat B$.  Since $A/\ker(\Ind\pi)$ is Morita equivalent to $B/\ker\pi$, Proposition~\ref{elementary} implies that $B/\ker\pi$ is elementary, and there is an isomorphism $\phi:B/\ker\pi\to \K(H)$.  The representation $\pi$ factors through a representation $\pi'$ of $B/\ker\pi$, and then $\pi'':=\pi'\circ\phi^{-1}$ is a representation of $\K(H)$ with $\pi(B)=\pi''(\K(H))$. Since every irreducible representation of $\K(H)$ is equivalent to the identity representation, we have
\[
\pi(B)=\pi''(\K(H))=\K(H_{\pi''})=\K(H_\pi),
\]
and $B$ is liminary. Symmetry gives the rest.
\end{proof}

One can quickly deduce from Proposition~\ref{prop-me-ccr} that $A$ is postliminary if and only if $B$ is, and that the Rieffel correspondence carries the largest liminary and postliminary ideals of $B$ to the corresponding ideals of~$A$. 

\begin{remark}
If one prefers to define postliminary algebras to be those for which $\pi(A)\supset\K(H_\pi)$ for every $\pi\in \hat A$ (it is actually a deep theorem that the two definitions are equivalent), then one can also prove directly that this property is preserved.  
\end{remark}

We learned the following lemma from Philip Green, and Zettl used a similar result in \cite{zettl2}.

\begin{lemma}\label{lem-green}
Let ${}_AX_B$ be  an imprimitivity bimodule and $\pi:B\to B(\H)$ a  representation of $B$.  For each $x\in X$, define $T_x=T_{x,\pi}:\H \to X\otimes_B \H$ by $T_x(h)=x\otimes h$ for $h\in H$.  Then $T_x^*(y\otimes h)=\pi(\langle x\,,\, y\rangle_B)h$ and
\[
T_x^*T_x=\pi(\langle x\,,\, x\rangle_B)\quad \text{and}\quad T_xT_x^*=\Ind\pi({}_A\langle x\,,\, x\rangle).
\]
\end{lemma}
\begin{proof}
For $x,y\in X$ and $h,k\in\H$ we have
\[
(T_x(k)\,|\, y\otimes h)=(\pi(\langle y\,,x\rangle_B)k\, |\, h)\\
=(k\, |\, \pi(\langle x\,,\, y\rangle_B)h),
\]
confirming the formula for $T_x^*$.  We have
$T_x^*T_x(h)=T_x^*(x\otimes h)=\pi(\langle x\,,\, x\rangle_B)h$.  Finally,
\begin{align*}
T_xT_x^*(y\otimes h)&=T_x(\pi(\langle x\,,\, y\rangle_B)h)=x\otimes \pi(\langle x\,,\, y\rangle_B)h\\
&=x\cdot \langle x\,,\, y\rangle_B\otimes h={}_A\langle x\,,\, x\rangle\cdot y\otimes h,
\end{align*}
which is by definition $\Ind\pi({}_A\langle x\,,\, x\rangle)(y\otimes h)$.
\end{proof}

\begin{cor}\label{cor-trace}
Let ${}_AX_B$ be  an imprimitivity bimodule and $\pi$ a  representation of $B$. For each $x\in X$,
\[
\tr(\pi(\langle x\,,\, x\rangle_B))=\tr(\Ind\pi({}_A\langle x\,,\, x\rangle)).
\]
\end{cor}

\begin{proof}
A slight modification of the proof that $\tr(T^*T)=\tr(TT^*)$ for $T\in B(\H)$ (for example, that given in \cite[Proposition~3.4.3]{ped-now}) shows that it holds also for $T\in B(H,K)$. Thus $\tr(T_x^*T_x)=\tr(T_xT_x^*)$, and the result follows from the Lemma.
\end{proof}

\begin{cor}\label{cor-wtf}
Let ${}_AX_B$ be  an imprimitivity bimodule, $\pi$ a  representation of $B$ and $x\in X$. Then  $\pi(\langle x\,,\, x\rangle_B)\neq 0$ if and only if $\Ind\pi({}_A\langle x\,,\, x\rangle)\neq 0$.
\end{cor}

Recall from \cite[4.5.2]{dix} that if $A$ is a $C^*$-algebra, then
\[
m(A):=\newspan\{a\in A^+:\pi\mapsto\tr(\pi(a))\text{\ is finite and continuous on $\hat A$}\}
\]
is an ideal, which we call the \emph{ideal of continuous-trace elements}.  Similarly,
\[
t(A):=\newspan\{a\in A^+:\pi\mapsto\tr(\pi(a))\text{\ is bounded on $\hat A$}\}
\]
is an ideal, which we call the \emph{ideal of bounded-trace elements} (\cite[\S 2]{mil}; see also \cite{ped69,per}).  The $C^*$-algebra $A$ has  \emph{continuous trace} if $\overline{m(A)}=A$, or \emph{bounded trace} if $\overline{t(A)}=A$.

\begin{prop}\label{btideal}
Let ${}_AX_B$ be  an imprimitivity bimodule.  Then the Rieffel correspondence carries $\overline{m(B)}$ to $\overline{m(A)}$ and $\overline{t(B)}$ to $\overline{t(A)}$.  In particular, $A$ has continuous trace if and only if $B$ has continuous trace, and $A$ has bounded trace if and only if $B$ has bounded trace.
\end{prop}

\begin{proof}
We prove the statement about $\overline{t(B)}$ and $\overline{t(A)}$; a very similar argument proves the analogous statement for $m(B)$.  Since every closed ideal $J$ in $A$ satisfies
\[
J=\textstyle{\bigcap}\{\ker\pi:\pi\in \hat A,\;\pi|_J=0\},
\]
it suffices to show that
\begin{equation}\label{eq-set}
\{\pi\in\hat B:\pi(\overline{t(B)})\neq\{0\}\}=\{\pi\in\hat B:\Ind\pi(\overline{t(A)})\neq\{0\}\}.
\end{equation}
Suppose $\pi(\overline{t(B)})\neq\{0\}$.  By \cite[Theorem~3.22]{tfb}, and then by polarisation,
\begin{align*}
\overline{t(B)}
&=\clsp\{\langle x\,,\, y\rangle_B:x,y\in X\cdot \overline{t(B)}\}\\
&=\clsp\{\langle x\,,\, x\rangle_B:x\in X\cdot \overline{t(B)}\}\\
&=\clsp\{\langle x\,,\, x\rangle_B:x\in X\cdot t(B)\}
\end{align*}
So there exists $x\in X\cdot t(B)$ such that $\pi(\langle x\,,\, x\rangle_B)\neq 0$.  Since $\langle x\,,\, x\rangle_B\in t(B)$, the function $\pi\mapsto \tr(\pi(\langle x\,,\, x\rangle_B)$ is bounded, and it follows from Corollary~\ref{cor-trace} that ${}_A\langle x\,,\, x\rangle\in t(A)$.  By Corollary~\ref{cor-wtf}, $\Ind\pi({}_A\langle x\,,\, x\rangle)\neq 0$. Thus $\Ind\pi(\overline{t(A)})\neq\{0\}$, and we have shown that the left-hand side of \eqref{eq-set} is contained in the right-hand side. A similar argument gives the other inclusion.
\end{proof}

It is well-known that a $C^*$-algebra need have no largest continuous-trace ideal: for example, in the algebra $A_3$ of \cite[Example~A.25]{tfb}, $\ker\pi_1$ and $\ker \pi_2$ are distinct maximal continuous-trace ideals whose intersection is $\overline{m(A_3)}$. The bounded-trace property is quite different, as the following result shows. It was first proved in \cite[Theorem~2.8]{ass}, but our argument seems more direct.

\begin{prop}
Every $C^*$-algebra $A$ has a largest bounded-trace ideal.
\end{prop}

\begin{proof}
We consider the set $\I$ of all closed ideals in $A$ which have bounded trace. Observe that if $a\in I^+$ belongs to $t(I)$, then $\pi(a)$ vanishes for $\pi\in\hat A\setminus \hat I$, and hence $a\in t(A)$ also. Let $J$ be the closure of $\newspan\bigcup_{I\in\I}I$. Then $J$ is an ideal in $A$, and 
$\newspan\bigcup_{I\in\I}t(I)$ is dense in $J$; since $\newspan\bigcup_{I\in\I}t(I)\subset t(J)$, $J$ has bounded trace.
\end{proof}

We can now deduce the following corollary from Proposition~\ref{btideal}.

\begin{cor}
Let ${}_AX_B$ be  an imprimitivity bimodule.  Then the Rieffel correspondence carries the largest bounded-trace ideal of $B$ to the largest bounded-trace ideal of $A$.
\end{cor}

\section{Multiplicity numbers}\label{secRob}

Our next goal is to prove that the upper and lower multiplicities of representations are preserved under Morita equivalence. We suppose that $B$ is a $C^*$-algebra, $\pi\in \hat B$ and $(\pi_\alpha)$ is a net in $\hat B$ such that $\pi$ is a cluster point of $(\pi_\alpha)$. We use the following characterisations of upper multiplicity relative to a net from \cite[Theorem~2.4]{ass}:
\begin{itemize}
\item[(rk)] $M_U(\pi,(\pi_\alpha))\leq k$ if and only if there exists $b\in B$ such that $\pi(b)\not=0$ and $\rank \pi_\alpha(b)\leq k$ eventually.

\item[(tr)] $M_U(\pi,(\pi_\alpha))\leq k$ if and only if there exists $b\in B^+$ such that $\pi(b)$ is a non-zero projection and $\tr \pi_\alpha(b)\leq k$ eventually.
\end{itemize}
Our statement of (rk) is sightly different from that of \cite[Theorem~2.4(iii)]{ass} in that we do not require $b$ to be positive, but the two are equivalent because the rank of $\pi_\alpha(b^*b)$ is the same as the rank of $\pi_\alpha(b)$.

\begin{theorem}\label{mupi}
Suppose ${}_AX_B$ is an imprimitivity bimodule. Let $\pi\in \hat B$ and let $(\pi_\alpha)$ be a net in $\hat B$ such that $\pi$  is a cluster point of $(\pi_\alpha)$. Then
\[
M_U(\pi, (\pi_\alpha))=M_U(X\dashind\pi,(X\dashind\pi_\alpha)).
\]
\end{theorem}

In the proof of the theorem we need the following standard lemma.

\begin{lemma}\label{Xn}
Suppose ${}_AX_B$ is an imprimitivity bimodule. The $n$-fold direct sum $X^n$ is an $M_n(A)$-$B$ imprimitivity bimodule with
\begin{align*}
(a_{ij})\cdot x&=\big(\textstyle{\sum_{j=1}^na_{ij}\cdot x_j}\big)_i,\\
x\cdot b&=(x_i\cdot b)_i,\\
{}_{M_n(A)}\langle x,y\rangle&=\big({}_A\langle x_i,y_j\rangle\big)_{i,j},\ \text{ and}\\
\langle x,y\rangle_B&=\textstyle{\sum_{i=1}^n\langle x_i,y_i\rangle_B},
\end{align*}
for $x=(x_i),y=(y_i)\in X^n$, $(a_{ij})\in M_n(A)$ and $b\in B$.
If $\pi:B\to B(H)$ is a representation, then there is a unitary isomorphism $U$ of $X^n\otimes_B H$ onto $(X\otimes_B H)^n$ such that $U(x\otimes h)=(x_i\otimes h)_i$, and $U$ intertwines $X^n\dashind \pi((c_{ij}))$ with the matrix $\big(X\dashind\pi(c_{ij})\big)_{i,j}$ in $M_n(B(X\otimes_B H))=B((X\otimes_B H)^n)$.
\end{lemma}

\begin{proof}[Proof of Theorem~\ref{mupi}]
It suffices to prove that
\begin{equation}\label{inequ}
M_U(X\dashind\pi,(X\dashind\pi_\alpha))\leq M_U(\pi, (\pi_\alpha));
\end{equation}
indeed, given \eqref{inequ}, we can apply it to the dual bimodule $\widetilde X$ to get
\begin{align*}
M_U(\pi,(\pi_\alpha))&=M_U(\widetilde X\dashind(X\dashind\pi),(\widetilde X\dashind(X\dashind\pi_\alpha)))\\
&\leq M_U(X\dashind\pi,(X\dashind\pi_\alpha)).
\end{align*}

First we suppose that $k:=M_U(\pi, (\pi_\alpha))$ is finite. By (tr), there exists $b\in B^+$ such that $\pi(b)$ is a non-zero projection and $\tr \pi_\alpha(b)\leq k$ eventually. Choose a continuous function $f\in C_c([0,\infty))$ such that $f(t)=0$ for $t$ near $0$, $f(1)=1$, and $f(t)\leq t$ for all $t\geq 0$. For large $\alpha$, $\pi_\alpha(b)$ is a positive compact operator (trace-class, in fact), and since $f(t)\leq t$ for all $t$ the spectral theorem implies that
\[
\tr(\pi_\alpha(f(b)))=\tr (f(\pi_\alpha(b)))\leq\tr(\pi_\alpha(b))\leq k.
\]
Since $f(0)=0$ and $f(1)=1$, we also have $\pi(f(b))=f(\pi(b))=\pi(b)$, so $\pi(f(b))$ is a non-zero projection.

The point of applying $f$ to $b$ is that $f(b)^{1/2}=f^{1/2}(b)$ lies in the Pedersen ideal $\kappa(B)$ of $B$, which  is contained in every other dense ideal of $B$ (see \cite[Theorem~5.6.1]{ped}). In particular, $\kappa(B)$ is contained in the ideal $\langle X,X\rangle_B$ spanned by the elements of the form $\langle x,y\rangle_B$, and thus there are finitely many elements $x_i, y_i\in X$ such that $f(b)^{1/2}=\sum_{i=1}^n\langle x_i,y_i\rangle_B$. Now
\begin{align*}
f(b)&=\Big(\sum_{i=1}^n\langle x_i,y_i\rangle_B\Big)\Big(\sum_{j=1}^n\langle x_j,y_j\rangle_B\Big)^*=\sum_{i,j=1}^n\langle x_i,{}_A\langle y_i,y_j\rangle\cdot x_j\rangle_B.
\end{align*}
The matrix $({}_A\langle y_i,y_j\rangle)_{i,j}$ is a positive element of the $C^*$-algebra $M_n(A)$ (see, for example, \cite[Lemma~2.65]{tfb}), and hence has the form $D^*D$ for some $D=(d_{ij})\in M_n(A)$. Thus
\begin{align*}
f(b)
&=\sum_{i=1}^n\Big\langle x_i,\sum_{k,j=1}^n d_{ki}^*d_{kj}\cdot x_j\Big\rangle_B\\
&=\sum_{k=1}^n\Big\langle \sum_{i=1}^n d_{ki}\cdot x_i,\sum_{j=1}^n d_{kj}\cdot x_j\Big\rangle_B\\
&=\sum_{k=1}^n \langle z_k,z_k\rangle_B,
\end{align*}
where we have written $z_k=\sum_{i=1}^n d_{ki}\cdot x_i$. Then with $z=(z_k)\in X^n$, we have realised $f(b)$ as the single inner product $\langle z,z\rangle_B$ for the $B$-valued inner product of Lemma~\ref{Xn}.

Recall that $\pi(\langle z,z\rangle_B)$ has the form $T_z^*T_z$ for the operator $T_z:h\mapsto z\otimes_B h$ of $H_\pi$ into $H_{X^n\dashind\pi}=X^n\otimes H_\pi$ (see Lemma~\ref{lem-green}). Since $\pi(\langle z,z\rangle_B)=\pi(f(b))=\pi(b)$ is a non-zero projection,
\[
T_zT_z^*=X^n\dashind\pi({}_{M_n(A)}\langle z,z\rangle)
\]
is also a non-zero projection. By Corollary~\ref{cor-trace}, for large $\alpha$ we have
\[
\tr\big(X^n\dashind\pi_\alpha({}_{M_n(A)}\langle z,z\rangle)\big)
=\tr\big(\pi_\alpha(\langle z,z\rangle_B)\big)\leq k.
\]
Thus (tr) gives $M_U(X^n\dashind\pi,(X^n\dashind\pi_\alpha))\leq k$.

To see that this statement passes to one about $M_U(X\dashind\pi,(X\dashind\pi_\alpha))$, we use (rk) to find $C=(c_{ij})\in M_n(A)^+$ such that $X^n\dashind\pi(C)\not=0$ and for large $\alpha$ we have $\rank(X^n\dashind\pi_\alpha(C))\leq k$. By Lemma~\ref{Xn}, $X^n\dashind\pi(C)$ is essentially the $n\times n$ matrix $\big(X\dashind\pi(c_{ij})\big)$, and we deduce that at least one entry $X\dashind\pi(c_{ij})$ in this matrix is non-zero. Since for large $\alpha$ we have
\[
\rank\big(X\dashind\pi_\alpha(c_{ij})\big)
=\rank \big(e_{ii}\big(X^n\dashind\pi_\alpha(C)\big)e_{jj}\big)
\leq \rank\big(X^n\dashind\pi_\alpha(C)\big)\leq k,
\]
we deduce from (rk) that $M_U(X\dashind\pi,(X\dashind\pi_\alpha))\leq k$, and we have proved \eqref{inequ} when $k=M_U(\pi, (\pi_\alpha))$ is finite.

As we commented earlier, this suffices to prove the theorem when $k$ is finite. In particular, if one of the upper multiplicities is finite, then so is the other; hence if one is infinite, the other must be too, and we also have equality when $M_U(\pi, (\pi_\alpha))$ is infinite.
\end{proof}

We now use Theorem~\ref{mupi} to obtain information about the lower multiplicity numbers $M_L$ defined in \cite[\S 2]{a} and \cite[\S 2]{as}. 

\begin{cor}
Suppose ${}_AX_B$ is an imprimitivity bimodule. Let $\pi\in \hat B$ and let $(\pi_\alpha)$ be a net in $\hat B$ such that $\pi$ is a cluster point of $(\pi_\alpha)$.  Then
\[
M_L(\pi,(\pi_\alpha))= M_L(\Ind\pi,(\Ind \pi_\alpha)).
\]
\end{cor}

\begin{proof}
Since $\Ind$ is a homeomorphism on spectra, $\Ind\pi$ is a cluster point of $(\Ind \pi_\alpha)$. So it suffices to show that $M_L(\Ind\pi,(\Ind \pi_\alpha))\leq M_L(\pi,(\pi_\alpha))$. By \cite[Proposition~2.3]{as} there exists a subnet $(\pi_{\alpha_i})$ of $(\pi_\alpha)$ such that $M_L(\pi,(\pi_\alpha))=M_U(\pi,(\pi_{\alpha_i}))$. Now Theorem~\ref{mupi} gives
\[
M_L(\Ind\pi,(\Ind \pi_\alpha))\leq M_U(\Ind\pi,(\Ind \pi_{\alpha_i}))=M_U(\pi,(\pi_{\alpha_i}))=M_L(\pi,(\pi_\alpha)),
\]
as required.
\end{proof}

\begin{cor}
Suppose ${}_AX_B$ is an imprimitivity bimodule and let $\pi\in \hat B$.
\begin{enumerate}
\item Suppose $\{\pi\}$ is not open in $\hat A$ (so that the lower multiplicity of $\pi$ is defined). Then $M_L(\pi)=M_L(\Ind\pi)$.
\smallskip
\item $M_U(\pi)=M_U(\Ind\pi)$.
\end{enumerate}
\end{cor}
\begin{proof}
It suffices to show that $M_*(\Ind\pi)\leq M_*(\pi)$.

(1) Since $\Ind$ is a homeomorphism, $\Ind \pi$ is not open in $\hat A$. By Proposi\-tions~2.2 and~2.3 of \cite{as}, there is a net $(\pi_\alpha)$ in $\hat B\setminus\{\pi\}$ converging to $\pi$ such that $M_L(\pi)=M_U(\pi,(\pi_\alpha))$. Now Theorem~\ref{mupi} gives
\[
M_L(\Ind\pi)\leq M_L(\Ind\pi,(\Ind\pi_\alpha))\leq M_U(\Ind\pi,(\Ind\pi_\alpha))\leq M_U(\pi,(\pi_\alpha))=M_L(\pi).
\]

(2)  By \cite[Proposition~2.2]{as} there exists a net $(\pi_\alpha)$ in $\hat B$ converging to $\pi$ such that
\[
M_U(\Ind\pi)=M_U(\Ind\pi,(\Ind\pi_\alpha)),
\]
which, by Theorem~\ref{mupi}, is $M_U(\pi,(\pi_\alpha))$, and thus less than or equal to $M_U(\pi)$.
\end{proof}

Recall that $A$ is a \emph{Fell algebra} if every $\pi\in\hat A$ is a \emph{Fell point}, that is, there exists $a\in A^+$  such that $\sigma(a)$ is a rank-one projection for all $\sigma$ near $\pi$ in $\hat A$.  It was observed in \cite[\S 3]{as2} that the Fell algebras are the algebras of type~$\text{I}_0$ studied in \cite[\S 6.1]{ped}. By \cite[Theorem~4.6]{a}, $A$ is a Fell algebra if and only if $M_U(\pi)=1$ for every $\pi\in\hat A$. Thus we have:

\begin{cor}
Suppose that ${}_AX_B$ is an imprimitivity bimodule.  Then $A$ is a Fell algebra if and only if $B$ is a Fell algebra.
\end{cor}

\section{Nuclearity}

Recall that a $C^*$-algebra $A$ is \emph{nuclear} if there is only one $C^*$-norm on the algebraic tensor product $A\odot C$ for every $C^*$-algebra $C$; the maximal tensor product $A\otimes_{\max}C$  the spatial tensor product $A\otimes_{\sigma}C$ (as defined and discussed in Appendix~B of \cite{tfb}, for example) then coincide.  Our goal in this section is to give a simpler and more direct proof of the following theorem of Zettl~\cite{zettl1} and Beer~\cite{beer}.

\begin{theorem}\label{thm-me-nuclear}
Suppose that ${}_AX_B$ is an imprimitivity bimodule.  Then $A$ is nuclear if and only if $B$ is nuclear.
\end{theorem}

For the proof of Theorem~\ref{thm-me-nuclear} we need the following lemma.

\begin{lemma}\label{lem-tfb+max}
Suppose that ${}_AX_C$  and ${}_BY_D$ are imprimitivity bimodules.  Then there are unique $(A\odot C)$- and $(B\odot D)$-valued inner products on the tensor product bimodule $Z=X\odot Y$ such that
\begin{gather}
\label{eq:418}
{}_{A\odot B}\langle x\otimes y,z\otimes w\rangle = {}_A \langle x,z\rangle\otimes {}_B\langle y,w\rangle \quad\text{and} \\
\langle x\otimes y,z\otimes w\rangle_{C\odot D} =\langle x,z\rangle_C \otimes\langle y,w\rangle_D,\label{eq:419}
\end{gather}
and  $Z$ is then both a pre-$(A\otimes_{\max}B)$--$(C\otimes_{\max}D)$ imprimitivity bimodule, and a pre-$(A\otimes_{\sigma}B)$--$(C\otimes_{\sigma}D)$ imprimitivity bimodule.
\end{lemma}

\begin{proof}
In \cite[Proposition~3.36]{tfb}, we show that \eqref{eq:418} and \eqref{eq:419} define positive sequilinear forms no matter what tensor product norm we use.  We also show that when we use the spatial norms, the module actions are bounded, so that $Z$ is a pre-$(A\otimes_{\sigma}B)$--$(C\otimes_{\sigma}D)$ imprimitivity bimodule.  To see that the same is true for the maximal norm, we only need to see that the module actions are bounded when the inner products are viewed as taking values in the maximal tensor products. We denote by $X\otimes_{\max}Y$ the Hilbert-module completion when $C\odot D$ has the maximal norm.

We begin by showing that the left action of $A$ is bounded. Consider the $A$--$M_n(C)$ imprimitivity bimodule $X^n$, defined as in Lemma~\ref{Xn} but with left and right swapped, and a typical element $\sum_{i=1}^nx_i\otimes y_i$ of $X\odot Y$. Since $A$ acts by bounded operators on $(X^n)_{M_n(C)}$ (see \cite[Lemma~3.7]{tfb}), there is a matrix $S=(s_{ij})$ in $M_n(C)$ such that
\begin{equation}\label{eq:420}
\big(\langle a\cdot x_{i},a\cdot x_{j}\rangle_C\big)_{i,j} = \|a\|^{2}\bigl(\langle x_{i},x_{j}\rangle_C\bigr)_{i,j} + S^{*}S.
\end{equation}
The matrix  $\bigl(\langle y_{i},y_{j}\rangle_D \bigr)$ is positive in $M_{n}(D)$, and therefore has the form $T^*T$ for some $T=(t_{ij})\in M_n(D)$. We now compute as follows:
\begin{align*}
\Big\langle\sum_{i} a\cdot x_{i}\otimes y_{i},&{} \sum_{i}a\cdot x_{i}\otimes y_{i} \Big\rangle_{C\odot D} = \sum_{i,j} \langle a\cdot x_{i},a\cdot x_{j}\rangle_C \otimes \langle y_{i},y_{j}\rangle_D \\
&= \sum_{i,j}\Bigl( \|a\|^{2} \langle x_{i},x_{j}\rangle_C +\sum_{k}s_{ki}^{*}s_{kj}\Bigr) \otimes \langle y_{i},y_{j}\rangle_D \\
&=\|a\|^{2} \Big\langle \sum_{i} x_{i}\otimes y_{i},
\sum_{j} x_{j}\otimes y_{j} \Big\rangle_{C\odot D} +\sum_{i,j,k,l} s_{ki}^{*}s_{kj}
\otimes t_{li}^{*}t_{lj} \\
&\le \|a\|^{2} \Big\langle\sum_{i} x_{i}\otimes y_{i},
\sum_{i} x_{i}\otimes y_{i} \Big\rangle_{C\odot D};
\end{align*}
since the term we threw away is positive in every $C^*$-completion of $C\odot D$, this last inequality holds in every completion, and in particular in the maximal tensor product.  A similar computation shows that $B$ acts by bounded operators on the second factor. The resulting homomorphisms of $A$ and $B$ into the $C^*$-algebra $\L(X\otimes_{\max}Y)$ have commuting ranges, and hence by \cite[Theorem~B.27]{tfb}  give a homomorphism of $A\otimes_{\max}B$ into $\L(X\otimes_{\max}Y)$, as required.
\end{proof}

\begin{proof}[Proof of Theorem~\ref{thm-me-nuclear}]
Suppose that $B$ is nuclear and $C$ is any $C^*$-algebra.  In view of \cite[Proposition~3.36]{tfb} and Lemma~\ref{lem-tfb+max}, the algebraic tensor product $X\odot C$ is both a
pre-$(A\otimes_{\max}C)$--$(B\otimes_{\max}C)$ imprimitivity bimodule, and a pre-$(A\otimes_{\sigma}C)$--$(B\otimes_{\sigma}C)$ imprimitivity bimodule.  In particular, the maximal norm of $t\in A\odot C$ is the operator norm of $t$ on $X\odot C$ viewed as a right Hilbert $(B\otimes_{\max}C)$-module, and the spatial norm is the operator norm of $t$ on $X\odot C$ viewed as a right Hilbert $(B\otimes_{\sigma}C)$-module. Since $B$ is nuclear, these norms coincide. Therefore, the maximal and spatial norms coincide on $A\odot C$.  Since $C$ is arbitrary, $A$ is nuclear.
\end{proof}

Presumably the following proposition is well-known, but we do not have a reference.

\begin{prop}
Every $C^*$-algebra has a largest nuclear ideal.
\end{prop}

For the proof, we need the following standard facts:
\begin{itemize}
\item[(a)] If $I$ is an ideal in $A$ and both $I$ and $A/I$ are nuclear, then so is $A$.
\smallskip
\item[(b)] If $A=\overline{\bigcup_i A_i}$ and each $A_i$ is a nuclear $C^*$-subalgebra of $A$, then $A$ is nuclear.
\end{itemize} 
The first of these is given a relatively elementary proof in \cite[\S B.53]{tfb}. For the second, let $B$ be a $C^*$-algebra and consider the canonical surjection $\phi:A\otimes_{\max}B\to A\otimes_{\sigma} B$. For each $i$ the norm inherited from $A\otimes_{\max}B$ is a $C^*$-norm on $A_i\odot B$, and hence the canonical map of $A_i\odot B$ into $A_i\otimes_\sigma B$ is isometric for this norm. Since the inclusion of $A_i\otimes_\sigma B$ in $A\otimes_{\sigma} B$ is isometric, it follows that $\phi$ is isometric on each subalgebra $\overline{A_i\odot B}$ of $A\otimes_{\max} B$, and hence $\phi$ itself is isometric.

\begin{proof}
We consider the collection $\mathcal{I}$ of nuclear ideals of $A$, which is nonempty because $\{0\}\in\mathcal{I}$. Property (b) above implies that chains in $\mathcal{I}$ have upper bounds in $\mathcal{I}$, and hence Zorn's lemma implies that $\mathcal{I}$ has maximal elements. If $I$ and $J$ are two maximal elements, then applying property (a) to the exact sequence
\[
0\to I\to I+J\to J/(I\cap J)\to 0
\]
shows that $I+J$ is nuclear, and hence maximality forces $I=J$. The unique maximal nuclear ideal is the one we want.
\end{proof}

Given the existence of the largest nuclear ideal, Theorem~\ref{thm-me-nuclear} immediately gives:

\begin{cor}
Suppose that ${}_AX_B$ is an imprimitivity bimodule. Then the Rieffel correspondence carries the largest nuclear ideal of $B$ into the largest nuclear ideal of~$B$.
\end{cor}

\end{document}